\documentclass[12pt,oneside]{amsart}
\usepackage{amscd,amsmath,amssymb,amsfonts,a4}
\usepackage{hyperref}
\usepackage{color}

\def\timesover#1#2#3{\ \xymatrix@1@=0pt@M=0pt{ _{#1}&\times&_{#2} \\& ^{#3}&}\ }
\def\otimesover#1#2#3{\ \xymatrix@1@=0pt@M=0pt{ _{#1}&\otimes&_{#2} \\& ^{#3}&}\ }

\usepackage[all]{xy}
\theoremstyle{plain}
\newtheorem{thm}{Theorem}
\newtheorem{lem}[thm]{Lemma}
\newtheorem{cor}[thm]{Corollary}
\newtheorem{prop}[thm]{Proposition}
\theoremstyle{definition}
\newtheorem{defn}[thm]{Definition}
\newtheorem{conj}[thm]{Conjecture}

\newtheorem{rmk}[thm]{Remark}
\newtheorem{rmks}[thm]{Remarks}

\newtheorem{eig}[thm]{Properties}

\numberwithin{thm}{section}
\numberwithin{equation}{section}

\newcommand{\ga}[2]{
\begin{gather}\label{#1}#2\end{gather} 
}

\newcommand{\surj}{\twoheadrightarrow}
\newcommand{\inj}{\hookrightarrow}

\newcommand{\Spec}{{\rm Spec \,}}

\newcommand{\sO}{{\mathcal O}}
\newcommand{\sP}{{\mathcal P}}

\newcommand{\sT}{{\mathcal T}}

\newcommand{\A}{{\mathbb A}}

\newcommand{\C}{{\mathbb C}}

\newcommand{\G}{{\mathbb G}}

\newcommand{\N}{{\mathbb N}}
\renewcommand{\P}{{\mathbb P}}
\newcommand{\Q}{{\mathbb Q}}

\newcommand{\Z}{{\mathbb Z}}

\newcommand{\id}{{\rm id}}
\newcommand{\FC}{\text{\sf FC}}
\def\tilde{\widetilde}

\begin{document}

\title[Packets]{Packets in Grothendieck's section conjecture}
\author{H\'el\`ene Esnault}
\address{
Universit\"at Duisburg-Essen, Mathematik, 45117 Essen, Germany}
\email{esnault@uni-due.de}
\author{Ph\`ung  H\^o Hai}
\address{
Universit\"at Duisburg-Essen, Mathematik, 45117 Essen, Germany
and Institute of Mathematics, Hanoi, Vietnam}
\email{hai.phung@uni-duisburg-essen.de}
\date{April 12, 2007}
\thanks{Partially supported by  the DFG Leibniz Preis and   the DFG Heisenberg program}

\parindent0cm
\parskip5pt

\begin{abstract}
The goal of this article is to show part of Grothendieck's section conjecture using the identification of sections with neutral fiber functors as defined in \cite{EP2}. 

\end{abstract}
\maketitle
\section{Introduction}

Let $U$ be an absolutely connected,  smooth scheme of finite type defined over field $k$ of characteristic 0. Then a section of the Galois group ${\rm Gal}(\bar k/k)$ into Grothendieck's fundamental group 
$\pi_1(U,\bar u)$ based at a geometric point $\bar u \to U$ is identified in \cite{EP2} with a neutral fiber functor $\rho$ of the Tannaka category ${\sf FC}(U)$ of finite connections. 
To a neutral fiber functor $\rho$, one defines  a $k$-form $s_\rho: U_\rho \to U$ of Grothendieck's
universal covering $U_{\bar u}\to U$ based a $\bar u$, which identifies  $\rho$ with the cohomological fiber functor $H^0(U_\rho, s_\rho^*(-))$. 

Grothendieck's section conjecture predicts a geometric description of sections ${\rm Gal}(\bar k/k)\to \pi_1(U,\bar u)$, under a geometric condition on $U$ and an arithmetic condtion on $k$. Without any condition, we can already say that the Tannaka description above is of geometric nature, if we think of neutral fiber functors as rational points of a gerbe. 
On the other hand, Tamagawa \cite{Ta} showed part of the conjecture following Grothendieck's suggestion. The aim of this article is twofold. We reprove the known part of Grothendieck's section conjecture using 
our identification of sections with neutral fiber functors as defined in \cite{EP2}. We then use our method to define the packets and to show the properties predicted by the conjecture.

A smooth absolutely connected curve over a field $k$ is called hyperbolic  if the degree of the sheaf of 1-differential forms with logarithmic poles at $\infty$ is strictly positive. Grothendieck's section conjecture \cite{GroFa} (see precise formulation in Conjecture \ref{conj:SC} of this article) predicts first that under the assumptions that $U$ be an hyperbolic curve
over $k$, a field of finite type over $\Q$, then  
rational points of $U$ inject into the set of sections. Call such sections geometric. 
Let $X(U_{\bar u})$ be the pro-system of smooth compactifications of the pro-system $U_{\bar u}$. Grothendieck predicts further that a geometric section yields  a unique ${\rm Gal}(\bar k/k)$-invariant point on $X(U_{\bar u})$, lying above a point in $U(k)$, and that a non-geometric section  also yields  a unique ${\rm Gal}(\bar k/k)$-invariant point on $X(U_{\bar u})$, thus lying above a point in $(X\setminus U)(k)$, where $X\supset U$ is the smooth compactification of $U$. 
Grothendieck mentions in \cite[p. 8]{GroFa} that those properties should be proven using a Mordell-Weil type argument. Indeed, it has been essentially worked out in \cite[Section~2]{Ta}. We show those properties anew in sections \ref{sec:i}, using our method. It becomes then a simple consequence of the definitions. 

According to Grothendieck, non-geometric sections should be  subdivided into packets $\sP_x, x\in (X\setminus U)(k)$ (\cite[p.9]{GroFa}). Each packet $\sP_x$ should have the cardinality of the continuum, and, for $x\neq y$, a section in $\sP_x$ should not be equivalent under conjugation with  $\pi_1(\bar U, \bar u)$ to a section in $\sP_y$. 

The Tannaka method we use, which relies on Nori's ideas \cite{N1}, \cite{N2},  on Deligne's non-neutral Tannaka theory \cite{DeGroth}, and on fiber functors at $\infty$ as defined in \cite[Section~15]{DeP} and \cite{KaGal}, allows us to define these packets and to show  the properties wanted.

In section \ref{sec:g}, we review Grothendieck's theory of the arithmetic  fundamental group
as used later on and  formulate precisely his  section conjecture. In section \ref{sec:s} we review our theory as developed in \cite{EP2} and suitable for our purpose here, and reformulate in this language Grothendieck's section conjecture. In section \ref{sec:h}, 
we provide a homological interpretation of the part concerning geometric sections. In section \ref{sec:i},  we show that
a section yields at most one  
 fixpoint on $X(U_{\bar u})$ and define the packets. In section \ref{sec:t}, we show that each packet has  the cardinality of the continuum. In fact, at each stage, we show what are the necessary conditions for those properties to hold. 

What remains to be understood is the more difficult part of Grothendieck's section conjecture, namely the existence of the fixpoint. Call it the existence conjecture. We make in section \ref{sec:e} a list of reductions. In particular, thinking of a fixpoint as a $k$-rational point of a $k$-form of $X(U_{\bar u})$, it would be enough to find a $K$-rational point for $K$ a function field appearing in the prosystem defining the $k$-form of $U_{\bar u}$ (Proposition \ref{Kpt=kpt}).  Furthermore, the existence conjecture suggests a weaker conjecture on extensions of fiber functors which, if possible to understand, would in general reduce the existence conjecture to an open in $\P^1$, and over a number field to $\P^1\setminus \{0,1,\infty\}$ (Proposition \ref{pe:8}).

{\it Acknowledgements:} 
Our first attempt was to try to contradict the existence conjecture using our method of neutral fiber functor \cite{EP2}. Pierre Deligne showed us our mistake in the computation of the fiber functor in the example we sent him. Thinking of his answer helped us to 
 develop  the formalism of this article. We thank him  for his interest and the clarity of his answer. We also thank the Deutsche Forschungsgemeinschaft (German Research Council) for its constant support.

\section{Grothendieck's section conjecture}\label{sec:g}
Let $k$ be a characteritic zero field, and fix an algebraic closure $\bar k$ of it together with the canonical inclusion $i:k\to \bar k$. Denote $\Gamma:=\text{Gal}(\bar k/k)$. Let  $\epsilon:X\to \Spec k$ be a absolutely connected
scheme $k$ and fix a closed geometric point $\bar x:\Spec  \bar k\to X$ of $X$.
Let $s_{\bar x}:X_{\bar x}\to X$ be the fundamental pro-\'etale covering of $X$, which is equipped with a $\bar k$-point $\tilde x$ with $s_{\bar x}(\tilde x)=\bar x$. Thus the arithmetic fundamental group $\pi_1(X,\bar x)$ of $X$, is, as a set, the fiber  $s_{\bar x}^{-1}(\bar x)$, in particular $\tilde x$ is identified with the unit element of $\pi_1(X,\bar x)$.

The morphism $\epsilon$ yields a homomorphism of fundamental groups
\ga{g:7}{\epsilon_*:\pi_1(X,\bar x)\to \pi_1(\Spec k,  i) = \text{Gal}(\bar k/k)=\Gamma.}
The assumption that $X$ is absolutely connected implies that $\epsilon_*$ is surjective.
The kernel of $\epsilon_*$ is isomorphic to $\pi_1(\bar X,\bar x)$, the geometric fundamental group of  $\bar X:=\Spec \bar k\times_{\Spec k}X$ with base point $\bar x$. 
That is, we have the following exact sequence
\ga{g:8}{1\to \pi_1(\bar X,\bar x)\to \pi_1(X,\bar x)\xrightarrow{\epsilon_*} \Gamma \to 1.}
In particular we have the following commutative diagram  
\ga{g:9}{\xymatrix{
&X_{\bar x}\ar[d]^{t_{\bar x},s_{\bar x}}\\
\Spec\bar k\ar[r]_{\id\times \bar x\qquad}\ar[ru]^{\tilde x}&\Spec\bar k\times X=\bar X.
}}
Here we omit the subindex $\Spec k$ for the fibre product of $k$-schemes over $\Spec k$.

Assume that $\bar x$ lies above a rational point $x:\Spec k\to X$, that 
is $\bar x=x\circ i$.
This yields a homomorphism $x_*:\pi_1(\Spec k,i)\to \pi_1(X,\bar x)$
which is a section of $\epsilon_*$ as a homomorphism of  pro-finite
groups. 

In general, let $y$ be an arbitrary $k$-rational point of $X$.
Then $\bar y:=y\circ i:\Spec \bar k\to X$ is a $\bar k$-point of $X$.
As above, $y_*$ is a section of $\epsilon_*:\pi_1(X,\bar y)\to \Gamma$.
Since $\pi_1(X,\bar x)\cong \pi_1(X,\bar y)$ by an isomorphism  which commutes with   $\epsilon_*$, we can consider $y_*$
as a section of $\epsilon:\pi_1(X,\bar x)\to \Gamma$, determined
up to an inner conjugation by an element of $\pi_1(\bar X,\bar x)$.

Grothendieck denotes the set of sections of $\epsilon_*:\pi_1(X,\bar x)\to \Gamma$
(as pro-finite group homomorphisms) up to conjugation by elements of
$\pi_1(\bar X,\bar x)=\text{Ker} \ \epsilon_*$
by $$\text{Hom-ext}_\Gamma(\Gamma,\pi_1(X,\bar x))$$  (``ext'' stands for
external; see  \cite[Eq.~(6)]{GroFa}). Notice that such a section  $\sigma:\Gamma\to \pi_1(X,\bar x)$
yields an action of $\Gamma$ on $X_{\bar x}$ making the morphism
$t_{\bar x}:X_{\bar x}\to \Spec\bar k$ $\Gamma$-equivariant.

\bigskip

\begin{defn} \label{defn:hyp}
Let $U$ be a smooth geometrically connected curve defined over a field $k$ of characteristic 0. Let $X\supset U$ be its smooth compactification. Then $U$ is said to be {\it elliptic} if its log-dualizing sheaf 
$\omega_{X}(\log (X\setminus U))$ has degree 0, and is said to be {\it hyperbolic} if $\omega_{X}(\log (X\setminus U))$ has strictly positive degree. Equivalently, $U$ is elliptic if the Euler characteristic 
$\chi(U):=2-g(X)-\#(X\setminus U)(\bar k)$ is $0$, and hyperbolic if it is strictly negative.

\end{defn}

Let  $k$ be a characteristic zero field,  and $U$ be a geometrically connected
curve over $k$.
Let $\bar u$ be
a $\bar k$-point of $U$. Then we have the exact sequence \eqref{g:8}  with $(X,\bar x)$ 
replaced by $(U,\bar u)$. Let $X(U_{\bar u})$  denote 
the compactification of $U_{\bar u}$ above the inclusion $U\inj X$, that is,
the pro-system of the compactifications of the finite \'etale coverings of
 $U$. It projects onto $X$, is \'etale above $U$ and ramifies along $(X\setminus U)$.

In a letter to Gerd Faltings dated June 27, 1983, Alexander Grothendieck proposed the following 
conjecture (see \cite[p.7-9]{GroFa}):
\begin{conj}[Grothendieck's section conjecture] \label{conj:SC}
Let $U$ be a hyperbolic curve defined over a field $k$ of finite type over $\Q$. Then 
\begin{itemize}
\item[(SC1)] The map 
\ga{}{U(k)\to \text{\rm Hom-ext}_\Gamma(\Gamma,\pi_1(U,\bar u)), \ y\mapsto y_*}
is injective. We call a section of the form $y_*$ a geometric
section.
\item[(SC2)] For a geometric section $y_*$, the resulting action of $\Gamma$
on $X(U_{\bar u})$ has a unique fixed point which lies above $y$.
\item[(SC3)] For a non-geometric section, the resulting action of 
$\Gamma$ on $X(U_{\bar u})$  has also a unique fixed point, which lies at 
infinity $X(U_{\bar u})\setminus U_{\bar u}$, and projects on a point in $(X\setminus U)(k)$.
This implies in particular that the set of non-geometric sections is divided
into  disjoint ``packets'' $\sP_x$, $x\in (X\setminus U)(k)$.
\item[(SC4)] Each packet $\sP_x$ has the cardinality of the continuum.
\end{itemize}
\end{conj}
Let us specialize to the case $X=U$. Then 
Grothendieck's section conjecture is equivalent to the following.
\begin{enumerate} \item[(SC)]
The map
\ga{}{X(k)\to \text{Hom-ext}_\Gamma(\Gamma,\pi_1(X,\bar x))}
is bijective. More precisely, let $\sigma$ be a section of $\epsilon_*:\pi_1(X,\bar x)\to \Gamma$. 
The scheme $X_{\bar x}$, considered as a $\Spec \bar k$-scheme by means of the morphism $t_{\bar x}$, has
a unique fixed point under the action of $\Gamma$ determined by means of $\sigma$.
Let $y$ be the image of this point in $X$. It is a $k$-rational point of $X$ and $\sigma=y_*$.
\end{enumerate}
\begin{rmk}\label{rmk:form}
We notice that the $\Gamma$-fixed points of $X_{\bar x}$ are nothing but
the $k$-rational points of the $k$-form of $X_{\bar x}$ determined
by the action of $\Gamma$. 
\end{rmk}

\section{The fundamental groupoid scheme}\label{sec:s}
As in the previous section, we fix a field $k$ of finite type over $\Q$
as well as an algebraic closure $\bar k$ and denote by $i:k\to \bar k$ the canonical inclusion. 
Let $X\xrightarrow{\epsilon}\Spec k$ be a smooth, geometrically connected scheme of finite type.
We denote by $\FC(X)$ the category of finite connections, defined in \cite[Section~2]{EP2}. This notion generalizes in characteristic 0 to smooth geometrically connected varieties  
Nori's notion
of finite bundles over proper geometrically
 connected schemes. As $X$ is geometrically connected, $\FC(X)$ is a Tannaka category over $k$ (as a fibre functor one can
always choose the tautological functor $\tau$ to ${\sf QCoh}(X)$ which assigns to a connection its underlying algebraic 
bundle).
Let $\rho:\FC(X)\to {\sf QCoh}(S)$ be a fibre functor, where $S$ is a $k$-scheme. Then Tannaka
duality yields a transitive $k$-groupoid scheme $\Pi_\rho\xrightarrow{(t,s)}S\times_kS$ \cite[Th\'eor\`eme~1.12]{DeGroth}.
Let $\Pi^\Delta_\rho$ denote the diagonal part of $\Pi_\rho$, which is a flat group scheme over $S$.
In \cite[Theorem~1.1]{EP2} we associate to $\rho$ a $\Pi^\Delta_\rho$-principal bundle $X_\rho\xrightarrow{(t_\rho,s_\rho)}
S\times_kX$. The original fibre functor $\rho$ is isomorphic to the cohomological functor 
$H_{DR}^0(X_\rho/S,s_\rho^*(-))$.
This construction is {\em functorial} in $\rho$, $S$ and $X$.
The following fibre functors are of special interest.

 Let us denote by ${\sf Fib}_k(X)$ the set of neutral fiber functors. 

Assume that $S=\Spec k$. Thus $\rho \in {\sf Fib}_k(X)$, hence $\Pi_\rho$ is a $k$-group scheme
and $s_\rho:X_\rho\to X$ is a principal under $\Pi_\rho$. In particular, let $x\in X(k)$ and $\rho:=\rho_x, \ \rho_x((V,\nabla))=V|_x$ be 
the fibre functor at $x$. Then our construction agrees with Nori's construction in case $X$ is proper. In this case, by definition $\Pi_{\rho_x}=s_\rho^{-1}(x)$ as a $k$-scheme.
The unit of $\Pi_{\rho_x}$ yields 
a $k$-rational point  $x_{\rho_x}$ of $X_{\rho_x}$ lying above $x$. This was a key point in Nori's theory. We use the following simplified notations:
\ga{s:0}{X_x:=X_{\rho_x}, \ \Pi(X,x):=\Pi_{\rho_x}, \ s_x:=s_{\rho_x}, \ \hat x:=x_{\rho_x}.}

For the tautological fibre functor $\tau$, we have
 $X_\tau=\Pi_\tau\xrightarrow{(t,s)}X\times X$. We call it the total fundamental groupoid scheme
of $X$.

\bigskip

Let $\bar x$ be a $\bar k$-rational point of $X$ and $\rho_{\bar x}$ be 
the fibre functor at $\bar x$. 
We call the $k$-groupoid scheme $\Pi_{\rho_{\bar x}}\xrightarrow{t,s}\Spec \bar k\times \Spec\bar k$
the fundamental groupoid scheme of $X$ with base point at $\bar x$ and denote it by $\Pi(X,\bar x)$.
By functoriality we have
\ga{s:1}{\xymatrix{
X_{\rho_{\bar x}}\ar[r]\ar[d]_{t_{\bar x},s_{\bar x}}\ar@{}[rd]|\Box&X_\tau\ar[d]^{t,s}\\
\Spec\bar k\times  X\ar[r]_{\quad \bar x\times\id}&X\times X
}\qquad
\xymatrix{\ar@{}[rd]|\Box
\Pi(X,\bar x)\ar[r]\ar[d]_{t,s}&X_\tau\ar[d]^{t,s}\\
\Spec\bar k\times \Spec \bar k\ar[r]_{\quad\bar x\times\bar x}&X\times X
}}
In particular $\Pi(X,\bar x)$ is the fibre of $X_{\rho_{\bar x}}\xrightarrow{s_{\bar x}}X$ at
$\bar x$. Thus we have the following commutative diagram
\ga{s:2}{\xymatrix{
&\ar[d]^{t,s} \Pi(X,\bar x)\ar[r]^{i}\ar@{}[rd]|{\Box}
& X_{\rho_{\bar x}}\ar[d]^{t_{\bar x},s_{\bar x}}\\
\Spec \bar k\ar[ru]_{e}\ar[r]_{\Delta\qquad}\ar@/^40pt/[rru]^{x_{\rho_{\bar x}}}&
\Spec\bar k\times\Spec\bar k\ar[r]_{\quad\id\times\bar x}&\Spec\bar k\times X 
}}
where $e$ denotes the unit of $\Pi(X,\bar x)$.
The composition $x_{\rho_{\bar x}}=i\circ e$ is a $\bar k$-point of $X_{\rho_{\bar x}}$ lying above $\bar x$.
It is shown in \cite[Theorem~1.2]{EP2} that $(X_{\rho_{\bar x}}, x_{\rho_{\bar x}})$ is the universal pro-\'etale covering of $X$
with respect to the fibre functor $F_{\bar x}$.  Furthermore, $\Pi(X,\bar x)$ considered as
a $\bar k$-scheme by means of the morphism $s:\Pi(X,\bar x)\to \Spec \bar k$, 
has the property that its $\bar k$-points build a profinite group which is identified with  $\pi_1(X,\bar x)$.

We shall therefore use the notation $X_{\bar x}: =X_{\rho_{\bar x}}$ 
and $\tilde x:=x_{\rho_{\bar x}}$.

Let $\bar y$ be another $\bar k$-point of $X$. Then $\Pi(X,\bar x)$ and $\Pi(X,\bar y)$ are isomorphic
as groupoid schemes acting upon $\Spec\bar k$.
Actually, any isomorphism is given by conjugating with
 an element of ${\rm Isom}^\otimes_{\bar k}(\rho_{\bar x},\rho_{\bar y})(\bar k)$,
which according to \cite[Theorem~3.2]{DeMil} is isomorphic (as a set) to $\Pi(X,\bar x)^\Delta(\bar k)$.

The following proposition was proved in \cite[Theorems~1.2, Theorem~1.3]{EP2}, except that we did not underline in loc.cit. that we keep the same $x$.
\begin{prop}\label{ps:1}
There is a one to one correspondence between conjugation classes of sections of $\epsilon_*:\pi_1(X,\bar x)\to \Gamma$
and equivalence classes of neutral fibre functors of ${\sf FC}(X)$. A section $\sigma$ of $\epsilon_*$ corresponds to a fibre functor $\rho$ if the $k$-form of $X_{\bar x}$ (see Remark \ref{rmk:form}) determined by $\sigma$ is isomorphic to $X_\rho$ over $X$. In particular a geometric section given by a point $x\in X(k)$ corresponds
to the fibre functor $\rho_x$ at $x$.
\end{prop}
\begin{proof}
According to the discussion above, $X_{\bar x}=X_{\rho_{\bar x}}$.
Given neutral fibre functor $\rho$ of ${\sf FC}(X)$,  then $\rho\otimes _k\bar k\cong
\rho_{\bar x}$, hence $X_\rho$ is a form of $X_{\rho_{\bar x}}$.
Conversely, each $k$-form of $X_{\rho_{\bar x}}$ determines a neutral fibre functor
by taking cohomology. On the other hand, each section of $\epsilon_*$ determines a
$k$-form of $X_{\bar x}$ and conversely.
\end{proof}
By means of this  proposition, Grothendieck's  section conjecture can be reformulated as follows.
\begin{conj}[Section conjecture reformulated]\label{cjs:2}
Let $U$ be a hyperbolic curve over a field $k$ of finite type over $\Q$ and $X$ be
its compactification.
 For each neutral fibre functor $\rho$
of ${\sf FC}(U)$ we denote by $X(U_\rho)$ the compactification of the $\Pi_\rho$-principal bundle
$s_\rho: U_\rho \to U$ with respect to $X$. Then
\begin{itemize}\item[(SC1)] The map $\delta:U(k)\to {\sf Fib}_k(X)$, $x\mapsto \rho_x$ is injective.
\item[(SC2)] For each $u\in U(k)$, the pro-curve $X(U_{\rho_u})$ has as unique $k$-rational point
which lies above $u$.
\item[(SC3)] Let $\rho\in{\sf Fib}_k(X)$ which is not geometric (i.e. not isomorphic to $\rho_u$ for any
$u\in U(k)$), then $X(U_\rho)$ has a unique $k$-rational point which lies above a point $x\in (X\setminus U)(k)$.
Denote the ``packet'' of such fibre functor in ${\sf Fib}_k(U)$ by $\sP_x$. Thus
\ga{}{{\sf Fib}_k(U)=\{\rho_u,u\in U(k)\}\bigsqcup \bigsqcup_{x\in (X\setminus U)(k)}\sP_x.}
\item[(SC4)] Each packet $\sP_x$ has the cardinality of the continuum.
\end{itemize}
\end{conj}
Let us specialize to the case $X=U$. Then 
Grothendieck's section conjecture reformulated is equivalent to the following.
\begin{enumerate} \item[(SC)]
The map
\ga{}{X(k)\to {\sf Fib}_k(X), \ x\mapsto \rho_x}
is bijective. 
\end{enumerate}

\section{A cohomology-like interpretation of (SC1) and (SC2)}\label{sec:h}
We continue to assume that $U$ is a curve over $k$, $k$ is of finite type over $\Q$.
Fix a neutral fibre functor $\rho$ for ${\sf FC}(U)$ (in particular assume the existence of such a functor).
Let $\eta$ be another neutral fibre functor for ${\sf FC}(U)$. Then according to \cite[Theorem~3.2]{DeMil}, the functor
$\text{Isom}^{\otimes}_k(\eta,\rho)$ is representable by a (profinite) scheme $E\to \Spec k$, which is a
principal bundle under both $\Pi_\rho$ (on the right) and $\Pi_\eta$ (on the left).
 Conversely, let $E\to \Spec k$ be a $\Pi_\rho$-principal bundle over
$\Spec k$, in particular $\Pi_\rho$ acts on $E$ from the right, hence $k[E]$ is
a representation of $\Pi_\rho$. Then we can define a fibre functor $\eta$ for ${\sf FC}(U)$ as follows:
\ga{h:1}{\eta(V):=(k[E]\otimes\rho(V))^{\Pi_\rho}, \text{ for all } V\in{\sf FC}(U)}
where $(-)^{\Pi_\rho}$ denotes the set of $\Pi_\rho$-invariants.
Then $E$ becomes a left $\Pi_\eta$-principal bundle. 

It is well known that $\Pi_\rho$-principal bundles over $\Spec k$ are classified by the pro-system of pointed sets 
\ga{h:11}{H^1_\text{\'et}(\Spec k,\Pi_\rho):=\varprojlim_{\Pi_\rho \surj G, G  \ {\rm finite}} H^1_{\text{\'et}}(k, G),} 
where
the distinguished element corresponds to a principal bundle $E$ with $E(k)\neq\emptyset$, thus $E\cong\Pi_\rho$ via the choice of the $k$-rational point. 
Consequently $\rho$ defines a bijective map $\delta_\rho:{\sf Fib}_k(U)\to H^1_\text{\'et}(\Spec k,\Pi_\rho).$

Now let us give ourselves $u\in U(k)$ and take $\rho=\rho_u$. We obtain a bijective map
\ga{h:2}{\delta_u:{\sf Fib}_k(U)\to H^1_\text{\'et}(\Spec k,\Pi(U,u)), \
v\mapsto E_v }
where $E_v$ is a short hand for the $k$-scheme representing the functor 
${\rm Isom}^{\otimes}(\rho_v, \rho_u)$. 

Observe that this map can be determined by the following diagram
\ga{h:3}{\xymatrix{E_v\ar[r]^{i_v}\ar[d]\ar@{}[rd]|\Box& X_{u}\ar[d]^{s_u}\\
\Spec k\ar[r]_v&X
}}

The geometric part of Grothendieck section conjecture (i.e. the part concerning geometric
sections) has the following cohomological interpretation.
Consider the ``sequence'' of $k$-schemes
\ga{h:4}{\Pi(U,u)\xrightarrow{i} U_{u}\xrightarrow{s_u} U,}
where $i$ is the closed embedding $s_\rho^{-1}(u)=\Pi(U,u)\subset U_u$.

Let $i^0$, $s_u^0$ denote the corresponding maps on the sets of $k$-points.
Define the ``cohomology set'' $H^1_\text{\'et}(k,U_{u})$
as the set of $k$-forms of $U_{u}\otimes \bar k$, and say that 
a $k$-form is trivial if it has a $k$-rational point. And define the map 
$$i^1:H^1_\text{\'et}(k,\Pi(U,u)) \xrightarrow{} H^1_\text{\'et}(k,U_{u})$$
that assigns to each $\Pi(U,u)$-principal bundle
$E\to \Spec k$ the $\Pi(U,u)$-principal bundle $U_E\to E$ defined as follows:
$E$ determines a fibre functor $\rho_E$ for ${\sf FC}(U)$, which in turn determines a $\Pi(U, \rho_E)$-principal bundle
 $U_{\rho_E} \to U$, which is a form of $U_{u}\otimes \bar k$. Then we have:
\begin{lem}\label{lh:1}
The sequence \eqref{h:4} induces the following ``long exact'' sequence of pointed sets
\ga{h:5}{e\to \Pi(U,u)(k)\xrightarrow{i^0}
 U_{u}(k)\xrightarrow{s_u^0} U(k)\xrightarrow{\delta_u}H^1_\text{\'et}(k,\Pi(U,u))
\xrightarrow{i^1} H^1_\text{\'et}(k,U_{u}).
}
\end{lem}
\begin{proof}
The exactness at $\Pi(U,u)(k)$ and at $U_{u}(k)$ is clear.
Let $v\in U(k)$ such that $\delta_u(v)$ is trivial in $H^1_\text{\'et}(k,\Pi(U,u))$,
which means that $E_v=\delta_u(v)$ has a $k$-rational point. According to
diagram \eqref{h:3}, $U_{u}$ has a $k$-rational point which lies above $v$.
Similarly, an element  $U_\rho$ of $H^1_\text{\'et}(k,\Pi(U,u))$ is mapped
to the trivial object of $H^1_\text{\'et}(k,X_{u})$ if and only if it has a $k$-rational point.
This mean $\rho$ is geometric, i.e. $U_\rho$ lies in the image of $\delta_u$.
\end{proof}
\begin{rmks}\label{rh:2} (a) The cohomological notations such as  $H^1_\text{\'et}(k,U_u)$ used in this section are just a support for thinking, there are not standard and also not properly defined as cohomologies. 

(b)
 The claim (SC1) of Grothendieck section conjecture is equivalent to saying that $\delta_u$ is injective,
or equivalently, that $s_u^0$ is trivial, i.e. maps $U_u(k)$ to the single point $u\in U(k)$,  or equivalently that $i^0$ is a bijective map.

(c) The claim (SC2) is equivalent to saying that $U_{u}(k)$ consists of a unique point above $u$.

(d)
If $U=X$ (i.e. $U$ is projective) then claim (SC3) is equivalent to saying that $\delta_u$
is surjective.
\end{rmks}

\section{The uniqueness of the fixed point and the packets}\label{sec:i}

In this section, using our method we provide new proofs of SC1 and SC2 in Conjecture \ref{cjs:2}, which were
 shown by Tamagawa in \cite[Section~2]{Ta}. 
This leads in a natural way  to the definition of the packets. 

\begin{thm} \label{thm:inj} Let $U$ be smooth geometrically connected curve which is elliptic or hyperbolic (see Definition \ref{defn:hyp}) and  defined over a field $k$ of finite type over $\Q$. 
Let $u\in U(k)$ be a $k$-rational point. Then the map 
\ga{i:1}{\delta_u:U(k)\to H^1_{\acute{e}t}(k,\Pi(U,u))}
is injective. 
\end{thm}
\begin{proof}
If $U$ is not projective, let $U\subset V$ be a larger open which is still
defined over $k$. 
It yields the restriction functor ${\sf FC}(V)\xrightarrow{{\rm rest}} {\sf FC}(U)$, which defines ${\sf FC}(V)$ as a full subcategory of ${\sf FC}(U)$. 
So as long as  $V$ still fulfills the assumption elliptic or hyperbolic, the theorem for $U$ follows from the theorem for $V$. So if the genus of the smooth compactification $X$ of $U$ is $\ge 1$, we can assume $U$ to be $X$, while if it is $0$, $U$ has to lie in a form of $\G_m$, which has to be $\G_m$ as we assume $U(k)\neq \emptyset$. In this case, we can assume $U$ to be $\G_m$. Let $J(V)$ be the Jacobian of $V$. So it is an abelian variety defined over $k$ in genus $\ge 1$ or the semiabelian variety $\G_m$ in genus $0$. We fix $u\in U(k)$ and denote by $j: U\to J(U)$ the cycle map sending $u$ to the unit $e$ of $J(V)$. To simplify notations, in the genus $0$ case, we assume $u=1$ so $j$ is then the identity. 
We denote by ${\sf FC}^\text{ab}(U)$ the full subcategory
of ${\sf FC}(U)$ consisting of those connections $V$  which split  into a direct sum
of rank 1 connections after base change by a finite extension of $k$. For a (neutral) fibre functor $\rho$, let
$\Pi^\text{ab}_\rho$ denote the corresponding Tannaka $k$-group scheme.  
It is the maximal commutative quotient $k$-group scheme of $\Pi_\rho$.

The morphism $j:U\to J(U)$ defines by pull-back the tensor functor $j^*:{\sf FC}(J(U))\to {\sf FC}(U)$,
which factors through the subcategory ${\sf FC}^\text{ab}(U)$. Therefore the resulting
homomorphism $j_*:\Pi(U,u)\to \Pi(J(U),e)$ factors through the homomorphism
\ga{i:2}{\bar j_*:\Pi^\text{ab}(U,u):=\Pi_{\rho_u}^\text{ab}\to \Pi(J(U),e).}
\begin{lem}\label{li:1}
The homomorphism $\bar j_*$ is an isomorphism.
\end{lem}
\begin{proof} Since ${\sf FC}(U)$ is compatible with base changes, so is ${\sf FC}^\text{ab}(U)$. Since $k$ has characteristic 0,  $\bar j_*$ is an isomorphism if and only if it is an isomorphism after $\otimes \bar k$. 
But over $\bar k$, $j$ induces an isomorphism $j^*: {\rm Pic}^0(J(U))\to {\rm Pic}^0(U)$ in genus $\ge 1$, while for $\G_m$, $j$ itself is the identity map. 
Therefore $\bar j^*$ yields an equivalence between ${\sf FC}(J(U))$ and ${\sf FC}^\text{ab}(U)$.
\end{proof}
\begin{prop}\label{pi:2}
Let $A/k$ be a semi-abelian variety defined over a field  $k$ of finite type over $\Q$. 
Then the map
\ga{i:4}{\delta_e: A(k)\to {\sf Fib}_k(A)}
is injective.
\end{prop}
\begin{proof}
For a natural number $m\in \N$, let us denote by $m \cdot A(k)$ the set of $k$-rational points of $A$ which are divisible by $m$ as $k$-rational points. 
Since $k$ is of finite type over $\Q$, there are no $k$-rational points which are infinitely divisible, thus 
\ga{i:5}{\bigcap_{m \in \N}  m\cdot A(k)=e.}
Let $s_e:A_e\to A$ be as usual 
the universal covering of $A$ associated to the fibre functor at $e$. By  Lemma \ref{lh:1}
the injectivity of $\delta_e$ is equivalent to the triviality of $s_e^0$, that is, all $k$-points of $A_e$
lies above $e$. Assume the contrary, that is $A_e$ has a $k$-point $b$ lying above $a\in A(k), \ a\neq e$. According to
\eqref{i:5}, there exists $m>0$ such that there is no $a'\in A(k)$ with $m\cdot a'=a$. Observe that
the covering $[m]:A\to A$ is a principal bundle under the $k$-group scheme  $A[m]:={\rm Ker}[m]$, and $e\in A(k)$ is a $k$-point lying above
$m\cdot e\in A(k)$ with respect to $[m]$. Therefore the universal covering
$s_e$ factors as
\ga{i:7}{\xymatrix{A_e\ar[r]^p\ar[rd]_{s_e}&A\ar[d]^{[m]}\\ & A
}}
Let $a':=p(b)\in A(k)$. Then $m\cdot a'=a$, a contradiction. Thus $s_e^0$ is the trivial map, hence $\delta_e$ is
injective.
\end{proof}
Proposition \ref{pi:2} finishes the proof of the Theorem.
\end{proof}
\begin{rmks} \label{rmk:injsharp}
Theorem \ref{thm:inj} shows SC1 in Conjecture \ref{cjs:2} under the weaker assumption on the geometry of $U$: it is enough for $U$ to be elliptic. The ellipticity assumption  is sharp. If $U=\A^1$, then certainly ${\sf FC}(\A^1)$ is the trivial category, thus all neutral fiber functors are equivalent to $M \mapsto H^0(\A^1, M)$, and $\delta_u$ can not be injective. 
It is also to be noted that for SC1, one only needs the maximal abelian subcategory of ${\sf FC}(U)$. 
\end{rmks}

\begin{thm}\label{tu:1} 
Let $U$ be smooth geometrically connected curve which is elliptic or hyperbolic (see Definition \ref{defn:hyp}) and  defined over a field $k$ of finite type over $\Q$. If 
$U$ is elliptic, we set $X=U$. If $U$ is hyperbolic, let $X$ be the smooth compactification of $U$.
Let $\rho$ be neutral fibre functor for ${\sf FC}(U)$, defining the principal $\Pi_\rho$-bundle $s_\rho: U_\rho\to U$ and the normalization $X(U_\rho)$ of $X$ in $k(U_\rho)$, defining $\bar{s}_\rho: X(U_\rho)\to X$ (see section \ref{sec:s}).   Then
\ga{u:1}{\#X(U_\rho)(k)\leq 1.}
\end{thm}
\begin{proof}
Assume the contrary that there exist $\alpha, \beta\in X(U_\rho)(k), \ \alpha\neq \beta$.
Let $a,b$ be respectively their images in $X(k)$. We use the notation $\rho_X$
 for the restriction of $\rho$ to ${\sf FC}(X)$.
Since ${\sf FC}(X)$ is a full subcategory of ${\sf FC}(U)$, we have the following
commutative diagram
\ga{u:2}{\xymatrix{
U_\rho\ar@{^(->}[r]\ar[d]_{s_\rho}& X(U_\rho)\ar@{->>}[r]\ar[d]_{\bar s_\rho}& X_{\rho_X}\ar[dl]^{s_{\rho_X}}\\
U\ar@{^(->}[r]&X
}}
Denote the images of $\alpha$ and $\beta$ in $X_\rho(k)$ by $ \tilde a$ and $\tilde b$. Their image in $X(k)$ are then 
$a$ and $b$. Thus we have
\ga{u:3}{\rho_X\cong \rho_a\cong\rho_b}
where $\rho_a$ and $\rho_b$ denote the geometric fibre functors at $a$ and $b$ for ${\sf FC}(X)$. Indeed, as recalled in section \ref{sec:s}, 
\ga{coh}{\rho_X  {\rm\ is \  equivalent\  to} \  {\sf FC}(X) \to {\sf Vec}_k, \ (V,\nabla)\mapsto H^0(X_{\rho_X}, s_{\rho_X}^*(V,\nabla))\\
{\rm thus \ to}  \ (V,\nabla)\mapsto  
V|_{\tilde{\alpha}}= V|_{a} \notag\\
{\rm  thus \ to} \  (V,\nabla)\mapsto  V|_{a}=V|_{b}.\notag}
We conclude by 
  Theorem \ref{thm:inj}  that $a=b$.

There exists a principal bundle $q: V\to U$ in the pro-system defining $U_\rho$ such that the images $a',b'$
of $\alpha$ and $\beta$ in the compatification $Y$ of $V$ are distinct. Thus
$\rho$ induces in a canonical way a fibre functor $\rho_V$ for ${\sf FC}(V)$. Indeed, $q_*N$ is an object of ${\sf FC}(U)$ whenever $N$ is an object of ${\sf FC}(V)$, and since $q^*q_*N\surj N$, the pull-back $N'$ of $N$ to $U_\rho$ is trivialized, and $\rho_V$ is defined by  $ N\mapsto H^0(U_\rho, N')$. It also shows that 
$V_{\rho_V}=U_\rho$. Consequently $X(U_\rho)$ is the compactification of $U_\rho$ with
respect to $V\inj Y$. On the other hand it is obvious that the pair $(V,Y)$ satisfies
the condition of the theorem.
Thus we can repeat the above argument to conclude that $a'=b'$. This is a contradiction.
This finishes the proof. 
\end{proof}

Theorem \ref{tu:1} suggests the following definition.
\begin{defn}\label{re:1} Let $U$ be a hyperbolic curve and $ X\supset U, k$ be as in 
Theorem \ref{tu:1}, where we allow $k$ to be any field of characteristic 0,  and let 
$x\in (X\setminus U)(k)$. We define the packet $\sP_x$ as the set of those 
fibre functors $\rho\in {\sf Fib}_k(U)$ which have the property that   the pro-scheme $X(U_\rho)$ has a $k$-point, and  this point lies above $x$.
\end{defn}
\begin{eig} \label{eig}
With the notations as in Definition \ref{re:1}, \eqref{coh} shows that if $\rho\in \sP_x$, then $\rho|_{{\sf FC}(X)}=\rho_x$. Vice-versa, if $\rho\in {\sf Fib}_k(U)$ such that $\rho|_{{\sf FC}(X)}=\rho_x$ for some $x\in (X\setminus U)(k)$, and $k$ is of finite type over $\Q$,  then to conclude that $\rho \in \sP_x$ would be a consequence of SC3 in Conjecture \ref{cjs:2}.
\end{eig}

\section{The tangential fibre functor}\label{sec:t}
  Let $U$ be a hyperbolic  curve (see Definition \ref{defn:hyp}) and  defined over a field $k$  of characteristic 0.
We keep the notation as in section \ref{sec:i}.

Let $x\in (X\setminus U)(k)$, with sheaf of maximal ideal $\frak{m}_x$.  We denote by $T_x=\Spec (\oplus_0^\infty \frak{m}_x^n/\frak{m}_x^{n+1})$ the fiber of the tangent bundle at $X$ in $x$, and by 
$T^0_x=T_x\setminus \{x\}=\Spec (\oplus_{-\infty}^{\infty}\frak{m}_x^n/\frak{m}_x^{n+1})$ 
the complement of the zero section. We denote by $K_x$ the local field at $x$, by $R_x$ its valuation ring, so after the choice of a local parameter $t$ at $x$, one has $K_x\cong k((t)), \ R_x\cong k[[t]]$. Denote
\ga{t:0}{S_x:=\Spec R_x,\quad S_x^0:= S_x\setminus \{x\}=\Spec K_x.}
Then there exists an exact tensor functor
\ga{tb}{{\sf FC}(S_x^0)\to {\sf FC}(T^0_x)}
constructed by Deligne \cite[Section~15]{DeP} and by Katz \cite{KaGal}, 
which is an equivalence of categories, and which inverts the natural restriction functor ${\sf FC}(T^0_x) \to {\sf FC}(S_x^0)$. 
So composing \eqref{tb} with the restriction functor $r_x:{\sf FC}(U)\to {\sf FC}(S_x^0)$ yields the tensor functor
\ga{ta}{DK_x: {\sf FC}(U)\to {\sf FC}(T_x^0).}
This allows   one to construct fibre functors for ${\sf FC}(U)$ by composing $DK_x$ with fibre functors
of ${\sf FC}(T_x^0)$. 
Since ${\sf FC}(S_x)$ is trivial, any fiber functor $\rho \in {\sf Fib}_k(U)$ obtained in this way has the property  
\ga{t1z}{\rho|_{{\sf FC}(X)}=\rho_x.}

Fibre functors of ${\sf FC}(T^0_x)$ can be  explicitly described. The choice of a local paramater $t$ at $x$ identifies $T_x^0$ with $\G_m= 
  \Spec k[t^{\pm 1}]$. Then ${\sf FC}(T_x^0)$ is
spanned (as an abelian category) by the connections
\ga{t:2}{L_a:=(\sO\cdot e_a,\nabla), \ \nabla(e_a)=\frac{adt}{t} e_a,\quad a\in \Q}
where $\sO:=\sO_{T^0_x}$. The tensor product for $L_a$ reads
\ga{t:3}{L_a\otimes L_b=L_{a+b}}
and there is an isomorphism
\ga{t:4}{\varphi_a:L_a\to L_{a+1}; \quad \varphi_a(e_a)=\frac{a+1}te_{a+1}.}
Thus ${\sf FC}(T_x^0)$ is indeed generated (as a tensor category) by $L_a, \ a\in \Q/\Z$.
The Tannaka group of ${\sf FC}(T^0_x)$ with respect to the fibre functor at $1\in \G_m$ identified with $T^0_x$ via $t$ is
\ga{t:5}{\Pi(T^0_x,1)\cong \varprojlim_n\mu_n=:\mu_\infty; \quad \mu_{m\cdot n}\xrightarrow{(-)^m} \mu_n.}
Consequently (see \eqref{h:11})
\ga{t:6}{{\sf Fib}_k(T^0_x)=H^1_\text{\'et}(k,\varprojlim_n\mu_n)=({\rm Kummer \ theory}) \varprojlim_n\frac{k^\times}{(k^\times)^n}.}
Notice that $k^\times$ is naturally embedded  in $\varprojlim_n\frac{k^\times}{(k^\times)^n}$ and its image
is precisely the set of geometric fibre functors of ${\sf FC}(T^0_x)$.  In particular we have
proved the following
\begin{cor}\label{ct:2}
The set ${\sf Fib}_k(T^0_x)$ has the cardinality of the continuum. 
\end{cor}

The universal covering $(T_x^0)_\rho$
associated to a fibre functor $\rho$ of ${\sf FC}(T_x^0)$ can be described as follows. Consider $\rho$ as
an element of $\varprojlim_n\frac{k^\times}{(k^\times)^n}$ and let $\rho_n\in k^\times$ be a representant
of the image of $\rho$ in $\frac{k^\times}{(k^\times)^n}$. Thus we have
\ga{t:7}{\rho_n\equiv \rho_{m\cdot n}\quad \mod (k^\times)^n.}
Then
\ga{t:8}{(T_x^0)_\rho=\Spec\big(\frac{k\left[t_{1/n}^{\pm 1}, n\in \N\right]}{(t_{1/m\cdot n})^m=
\frac{\rho_{m\cdot n}}{\rho_n}t_{1/n}}\big).}
The projection $(T^0_x)_\rho\to T^0_x$ is defined  by $t\mapsto t_1$.

\begin{rmk} \label{rmk:exsharp} This remark echoes Remarks \ref{rmk:injsharp}. 
The normalization $T_x(T_x^0)$ of $T_x$ in $k((T_x^0)_\rho)$  is
\ga{t:9}{T_x((T^0_x)_\rho)=\Spec\big( \frac{k\left[t_{1/n},n\in \N\right]}{(t_{1/m\cdot n})^m=
\frac{\rho_{m\cdot n}}{\rho_n}t_{1/n}}).}
It has a unique $k$-point (which lies above $0\in T_x$) if and only if $\rho$ is not geometric.
We see in this way that the condition on $U$ being elliptic in Theorem \ref{tu:1} is sharp. 
\end{rmk}

We denote by $\sT_x$ the set of fibre functors of ${\sf FC}(U)$ obtained by composing $DK_x$ with 
a fibre functor $\rho$ of ${\sf Fib}_k(T_x^0)$, and use the  notation 
\ga{t:10}{\tau_{x,\rho}:=\rho\circ DK_x.}

Let us choose a $k$-rational point of $T_x^0$, which we denote by $1$, as it is the same as choosing a parameter $t\in \frak{m}_x/\frak{m}_x^2\setminus \{0\}$ at $x$ and the point defined by $t=1$.  
Tannaka duality for
\ga{t11}{\xymatrix{{\sf FC}(U)\ar[r]^{DK_x}\ar[rd]_{\tau_{x,\rho_1}}& {\sf FC}(T_x^0)\ar[d]^{\rho_1}\\
&{\sf Vect}_k 
}}
where $\rho_1$ denotes the fibre functor at $1\in T_x^0$, yields a group scheme homomorphism
\ga{t:12}{DK_x^*:\Pi(T^0_x,1)\to \Pi(U,\tau_{x,\rho_1}).}
\begin{lem}\label{lt:3}
The homomorphism $\text{DK}_x^*$ in \eqref{t:12} is injective.
Consequently each object of ${\sf FC}(T^0_x)$ is isomorphic to a direct summand
of the image under $\text{DK}_x$ of an object of ${\sf FC}(U)$.
\end{lem}
\begin{proof}
The homomorphism ${DK}_x^*$ respects base change. 
So we may assume that  $k$ is of finite type over $\Q$, thus in particular, it is  embeddable  in $\C$. 
Then \eqref{t:12}  is compatible
with the homomorphism of topological fundamental groups
\ga{t:13}{\pi_1(\C^\times,1)\to \pi_1(U(\C), t)}
after the choice of an embedding $k\to \C$. Here $\pi_1(U(\C), t)$ is the topological fundamental group based at the tangent vector $t$ at $x$ in the sense of Deligne (see
\cite[Section~15]{DeP}). 
Ellipticity or hyperbolicity  of $U$ implies that  \eqref{t:13} is injective. 
Since by uniformization theory, the topological fundamental groups of $U(\C)$ and of $\C^\times$ lie inside the complex points of an algebraic group, thus are residually finite,  $DK^*_x$ is injective as well. 

The last claim follows from \cite[Theorem~2.11]{DeMil} and the fact
that ${\sf FC}(U)$ is a semi-simple abelian category.\end{proof}
We mention the following result of \cite[Theorem~5.7]{EP2}.
\begin{prop}\label{pt:3}
Let $\tau:=\tau_{x,\rho}\in \sT_x$. Then $X(U_\tau)$
has a $k$-point lying above $x$. Consequently $\sT_x\subset \sP_x$.
\end{prop}

The aim of the rest of this section is to show 
\begin{thm} \label{thm:T=P}
Let $U$ be a hyperbolic curve (see Definition \ref{defn:hyp}) defined over a field $k$ of characteristic $0$. Then $\sT_x=\sP_x$ for any $x\in (X\setminus U)(k)$. 
\end{thm}
We first prove intermediate statements. 

 Fix $\rho\in{\sf Fib}_k(U)$. Then each finite full tensor subcategory $S\subset {\sf FC}(U)$ defines a finite quotient $\Pi_\rho\surj G:=\Pi_{\rho,S}$,
 and a principal bundle $U_G:=U_{\rho,S}\xrightarrow{p} U$ under group $G$.
 Let $X(U_G)$ be the compactification of $U_G$, which then projects on $X$.

We set
$$S_{G,x}:=S_x\times_XX(U_G)\ \text{ and }\ S_{G,x}^{0}:=S_x^0\times_UU_G $$
where $S_x, S^0_x$ are defined in \eqref{t:0}.

 \begin{prop}\label{pt:5}
The (fixpoint free) action of $G$ on $U$
extends to a (not necessarily fixpoint free) action
\ga{t:13a}{
\xymatrix{\ar[dr]_{p_1} G\times X(U_G)\ar[r]^{\bar \mu} &   X(U_G) \ar[d]^{\bar p}\\
& X  } \xleftarrow{{\rm rest}} 
\xymatrix{\ar[dr]_{p_1} G\times S_{G,x} \ar[r]^{\bar \mu} &  S_{G,x} \ar[d]^{\bar p}\\
& S_x  }\xleftarrow{{\rm rest}}
\xymatrix{\ar[dr]_{p_1} G\times S_{G,x}^0 \ar[r]^{\mu} &  S_{G,x}^0 \ar[d]^{ p}\\
& S^0_x  }}

\end{prop}
 \begin{proof} Denote for simplicity $V:=U_G$  and $Y:=X(U_G)$.	
We first start with the action on $Y$. The composition of the two maps 
$G\times V\xrightarrow{\mu} V \to Y$ factors (uniquely) through 
$G\times Y$ as $G\times Y$ is a smooth curve and $Y$ is projective
 (it is a removable singularity). The same argument yields that the two maps 
$G\times G\times V\rightrightarrows G\times V\to G\times Y$
defined by $\mu(gh,v)$ and $ \mu(g,  \mu(h,v))$ factor uniquely  
through $G\times G \times Y\rightrightarrows Y$. Since $\mu(gh,v)= \mu(g,  \mu(h,v))$ on $V$, equality holds on $Y$ as well. This shows the existence of the extension of $\mu$ to $\bar \mu$ on the left triangle. On the other hand, since $p_1= p\circ \mu$ on $V$, one has $p_1=\bar p\circ \bar \mu$ on $Y$. This shows the left triangle. It also shows that the fiber $X(U_G)_x=\bar p^{-1}(x)$ is left invariant by $\bar \mu$. This finishes the proof. 
 \end{proof}

Let $\rho\in\sP_x$ and $\alpha\in X(U_\rho)$ be the unique $k$-point,
which lies above $x$.
Recall from Properties \ref{eig} that
$\rho|_{{\sf FC}(X)}=\rho_x$ and that ${\sf FC}(X)$ is a full tensor subcategory
of ${\sf FC}(U)$. Thus we have a surjective homomorphism $\pi(U, \rho)\to \pi(X, \rho_x)$,
the  kernel of which will be denoted by $K$.

\begin{prop}\label{tt:factor} The map $\bar s_\rho:X(U_\rho)\to X$ factors as
\ga{t:25}{\xymatrix{\ar@/_2pc/[dd]_{\bar{s}_\rho} X(U_\rho) \ar[d]^{\bar q_\rho} & \alpha \ar[d]\\
 X_{\rho_x} \ar[d]^{s_{\rho_x}} & x_{\rho_x} \ar[d] \\
X & x
}}
with the property
\ga{t:26}{\bar q_\rho^{-1}(x_{\rho_x})=\alpha}
  i.e. $ \bar q_\rho$  fully  ramifies  at $\alpha$. 
\end{prop}
\begin{proof}
 Fixing a finite full tensor subcategory $S$ of
${\sf FC}(U)$ as in the discussion preceeding Proposition \ref{pt:5} defines
 a finite quotient group scheme $G$ of $\pi(U, \rho)$. Recall that
 $K:=\text{Ker}(\pi(U,\rho)\to \pi(x,\rho_x)$.
Denote by $A$ the image of $K$ in $G$, and consider the exact sequence
\ga{t:17}{1\to A\to G\to H\to 1}
The group $H$ is isomorphic to the Tannaka group
of the restriction of $\rho$ to the intersection of $S$ with ${\sf FC}(X)$.
Thus $X_H$ is the compactification $X(U_H)$ of $U_H$  and
we have the following commutative diagram
\ga{t:18}{\xymatrix{\ar[d]_{A}^q U_G \ar[r] & X(U_G) \ar[d]^{\bar q} & x_G \ar[d]\\
\ar[d]_{H} U_H \ar[r] & X_H=X(U_H) \ar[d] & x_H \ar[d]\\
U \ar[r] & X & x
}}
where $q:U_G\to U_H$ is an $A$-principal bundle and $x_G, x_H$ are images of
the point $\alpha \in X(U_\rho)(k)$ in $X(U_G)$ and $X(U_H)$, respectively.

Proposition \ref{pt:5} applied to the principal bundle $U_G\xrightarrow{q}U_H$
yields an action
\ga{t:19}{A\times \bar q^{-1}(x_H)\to \bar q^{-1}(x_H) }
of the $k$-algebraic group scheme $A$ on the $k$-scheme $\bar q^{\ -1} (x_H)$. 

One has from \eqref{t:18} that $x_G\in \bar q ^{-1}(x_H)$ and that
$H$ acts freely on $X_H$. Consequently $B:={\rm Stab}(x_G)\subset G$
projects trivially onto $H$, thus $B\subset A$. This yields the exact sequence 
\ga{t:20}{
1\to A/B \to G/B\to H\to 1}
together with a  principal bundle
\ga{t:21}{A/B\times \bar q^{-1}(x_H)\to x_H.}
Thus one has a factorization
\ga{t:22}{\xymatrix{\ar[d]_B^{q''} U_G \ar[r] & X(U_G) \ar[d]^{\bar q''} & x_G \ar[d]\\
\ar[d]_{A/B}^{q'} U_{G/B} \ar[r] & X(U_{G/B})\ar[d]^{\bar q'} & x_{G/B} \ar[d]\\
U_H \ar[r] & X_H & x_H
}
}
Denoting by $e$ the index of ramification of $x_G\otimes \bar k$, one has that $|B(\bar k)|=e$ while $ |(A/B)(\bar k)|=| \bar q^{-1}(x_H)(\bar k)|$. We conclude that 
\ga{t:23}{|(\bar q')^{-1}(x_H)(\bar k)|= | \bar q^{-1}(x_H)(\bar k)|=|(A/B)(\bar k)|} 
and $\bar q'$ is \'etale. Thus $A/B$ factors as a quotient of $\pi(X,\rho_x)$. 
By definition \eqref{t:17}, one concludes $A/B=1$. This shows 
\ga{t:24}{\bar q^{-1}(x_H)=\{x_G\}.}
Passing to the pro-system we deduce the claims of the proposition.
  \end{proof}
\begin{thm}\label{tt:6}
With the assumption of Theorem \ref{thm:T=P}, for any $\rho \in \sP_x$ there is a factorization 
\ga{t:14}{\xymatrix{\ar[d]_{\rho} {\sf FC}(U) \ar[r]^{r_x} & {\sf FC}(S_x^0) \ar[dl]^{\exists \eta}\\
 {\sf Vec}_k
}
}
\end{thm}
\begin{proof}
We adopt the notation of Proposition \ref{tt:factor}.
 By Hensel's lemma, one has a unique lifting
\ga{t:27}{\xymatrix{ & X_{\rho_x}\ar[d]^{s_{\rho_x}}\\
\ar[ur]^i S_x \ar[r] & X}
}
 with $x_{\rho_x}\in i(S_x)$. 
We define
\ga{t:28}{\sigma: \frak{S}_x^0:=S_x^0\times_{i(S_x)} U_\rho\to S_x^0\\
\notag
\xymatrix{\frak S_x^0\ar[r]\ar[d]_\sigma\ar@{}[rd]|{\Box}&U_\rho\ar[d]^{\bar q_\rho}
\ar[rd]^{s_\rho}\\
S_x^0\ar[r]_i& X_{\rho_x}\ar[r]_{s_{\rho_x}}&X
}}

According to Lemma \ref{lt:3}, every object $N$ in ${\sf FC}(S_x^0)$ is a summand of 
the restriction to $S_x^0$ of an object $M$ of ${\sf FC}(U)$. Since $s_\rho^*(M)$ is trivializable, 
so is $\sigma^{*}(N)$. On the other hand, \eqref{t:26} implies 
\ga{t:29}{H^0(\frak{S}_x^0)=k.}
Consequently, $\rho$ extends to ${\sf FC}(S_x^0)$ by setting
\ga{t:30}{\rho(N)=H^0(\sigma^*(N)).}
\end{proof}
\begin{proof}[Proof of Theorem \ref{thm:T=P}] 
By means of Theorem \ref{tt:6} and of the equivalence
in \eqref{tb}, the claim of Theorem \ref{thm:T=P} amounts
to the uniqueness of the functor $\eta$ in diagram \eqref{t:14} for any
given functor $\rho\in\sP_x$. Thus, let $\rho$ and $\eta$ be as in \eqref{t:14}
and let $\varphi$ be another fibre functor in ${\sf Fib}_k(S_x^0)$,
such that $\rho\cong \varphi\circ r_x$.
Via  the equivalence ${\sf FC}(S_x^0)\equiv {\sf FC}(T_x^0)$, we can consider $\varphi$ as a functor
in ${\sf Fib}_k(T_x^0)$. Now, as in the proof of Proposition \ref{pt:3}, i.e. the proof of
\cite[Theorem~5.7]{EP2}, there exists a map
\ga{}{f:(S_x^0)_{\varphi}\to  U_\rho}
that determines the $k$-point of $X(U_\rho)$, see diagram (5.11) of [loc.cit].
The uniqueness of the $k$-point of $X(U_\rho)$ shows that $f$ is compatible
with the morphisms $i$ of \eqref{t:27} and $\bar q_\rho$ of \eqref{t:25} in the sense
that we have the commutative diagram
\ga{}{\xymatrix{
(S_x^0)_{\varphi} \ar[r]^f\ar[d]_{s_\varphi} & U_\rho \ar[d]^{\bar q_\rho}\\
S_x^0\ar[r]_i &X_{\rho_x}
}}
It follows from the universal property of $\frak S_x^0$, i.e. \eqref{t:28}, that there
exists a map
\ga{}{\theta:(S_x^0)_{\varphi}\to\frak S_x^0}
which is compatible with the maps to $S_x^0$ and $U_\rho$.
That is, there exists a natural tensor transformation $\varphi\to \eta$,
which is then automatically a natural isomorphism (cf. \cite[Proposition~1.13]{DeMil}).
Thus $\eta$ and $\varphi$ are isomorphic.
\end{proof}
\begin{cor}\label{ct:cardinality}
Let $\eta$ and $\eta'$ be non-isomorphic fibre functors of ${\sf FC}(T_x^0)$.
Then $\eta\circ DK_x$ and $\eta'\circ DK_x$ are not isomorphic.
Consequently, if $k$ is of finite type over $\Q$ then $\sP_x$ has the cardinality of the continuum.
\end{cor}

\section{Remarks towards the existence}\label{sec:e}
In this section, we gather a few facts  towards SC3 in Conjecture \ref{cjs:2}.

\begin{prop} \label{Kpt=kpt}
Let $U$ be a geometrically connected curve defined over a characteristic 0 field $k$.  Let $\rho \in {\sf Fib}_k(U)$ and $X(U_\rho)$ be as in Theorem \ref{tu:1}.  Then, if $U$ is hyperbolic (see Definition \ref{defn:hyp}), then  $X(U_\rho)(k)=X(U_\rho)(L)$ for $L=k(V)$ where $V\to U$ is any member of the pro-system $s_\rho: U_\rho \to U$.

\end{prop}
\begin{proof}
Let $a: \Spec L\to X(U_\rho)$, Then for $W\to U$ in the pro-system defining $U_\rho$, with $W\neq V$, the hyperbolicity condition implies that the genus of the compactification of $W$ is strictly larger than the genus of the compactification of $V$. Therefore the induced point $\Spec L\to X(W)$ has to be a closed point. Let $\kappa(a)$ be its residue field, thus $k\subset \kappa(a)$ is finite, and $\kappa(a) \subset L$.  By the connectivity of $V$, this is only possible if $\kappa(a)=k$. This finishes the proof. 
\end{proof}

The following remark is contained in \cite[Corollary~2.10]{Ta}, and seen here quite directly with our method. 
\begin{rmk}\label{pe:1}
Let $U$ be a hyperbolic curve  defined over a field $k$ of finite type over $\Q$. 
 Let $\rho$ be a neutral fibre functor of ${\sf FC}(U)$.
Then $X(U_\rho)$ has a $k$-point if and only if for any principal bundle $V$ in the pro-system
defining $U_\rho$, $X(V)(k)\neq \emptyset$, where $X(V)$ is the compactification of $V$.
\end{rmk}
\begin{proof}
The ``only if'' part is obvious since the image of the $k$-point in $X(U_\rho)$ on each $X(V)$
is a $k$-point of $X(V)$. We prove the ``if'' part.

Assume that $X$ has genus $\geq 2$. Then, as a consequence of Faltings' theorem asserting that a genus $\ge 2$ curve over a finite field has at most finitely many rational points,   there is a
Hausdorff topology on $k$ that induces a topology on $X(k)$ making it
a compact space. 
For each $S\subset {\sf FC(U)}$ finite full subcategory we denote by $U_S$ the
principal bundle constructed from $S$ and $\rho$.
 We call a point $a\in X(k)$ $S$-good
if  $a$ lies in the image of $X(U_S)(k)$, which, by assumption, is not empty.
By assumption on the genus of $X$, the genus of $X(U_S)$ is also $\geq 2$.
Thus for each $S$, the set of $S$-good points is not empty, and is compact,
being the image of a compact set. For $S$, $S'$ finite 
subcategories, the set of points which are both $S$-good and $S'$-good is not empty neither: just
consider $S\cup S'$.
Thus the set of $\rho$-good points in $X(k)$, i.e. 
those which are $S$-good for any $S$,
is not empty. This can be repeated for any principal bundle $V$ of $U$ with the fibre functor induced 
from $\rho$ (see proof of Theorem \ref{tu:1}). 
Now we can conclude the existence of the $k$-point in the universal covering 
$X_\rho$ as points in the limit of the pro-system of $\rho$-good points in each covering $Y$ 
of $X$.

If $X$ has genus 0 then we can still find in the pro-system defining $U_\rho$ a principal bundle 
$V$ with compactification $Y$ of genus at least 2. We can then replace $(U,X)$ by $(V,Y)$ without
lost of generality.
\end{proof}

\begin{cor}\label{ce:3} 
Claim SC3 of Conjecture \ref{cjs:2}  holds for all hyperbolic curves $U$ over $k$ if and only if it does
 for the subfamily of all $U$ with $X(k)=\emptyset$, that is in this case
${\sf Fib}_k(U)=\emptyset$.
\end{cor}
\begin{proof}
Let $U$ be given. Assume $\rho$ given and $X(U_\rho)(k)=\emptyset$. Then, according to Remark \ref{pe:1}
 there is a $V\to U$ in the pro-system, with compactification $Y\to X$, 
such that $Y(k)=\emptyset$, as $Y$ is in the pro-system of $X(U_\rho)$. So
$V$ violates the conjecture. 
\end{proof}
\begin{prop}\label{pe:2} The assumptions are as in SC3 of Conjecture \ref{cjs:2}. Then 
$|X(U_\rho)(k)|= 1$
 and only if $|X_L((U_L)_{\rho_L})(L)|= 1$  for some $L\supset k$ finite.
 Consequently if the section conjecture holds for $U_L$ then it holds for $U$.
\end{prop}
\begin{proof}
Assume  $|X_L((U_L)_{\rho_L})(L)|= 1$. One has a fiber square
\ga{2}{\xymatrix{  \ar @{} [dr] |{\square}\ar[d] (U_L)_{\rho_L} \ar[r] & U_\rho \ar[d]\\
U_L \ar[r] & U.
}
}
Since the normalization of $X$ in $k(U_\rho)\otimes_k L$ is the normalization of $X_L$ in $k(U_\rho)\otimes_k L$, \eqref{2} induces the fiber square 
\ga{3}{\xymatrix{  \ar @{} [dr] |{\square}\ar[d] X_L((U_L)_{\rho_L}) \ar[r] & X(U_\rho) \ar[d]\\
X_L \ar[r] & X.
}
}
Assume $X_L((U_L)_{\rho_L})(L)=\{\alpha\}$. Then  $X((U)_{\rho})(L)=\{\alpha\}$. 
If there was no $a\in X(U_\rho)(k)$ below $\alpha$, 
then $X(U_\rho)_L(L)$
 would contain the conjugates of $\alpha$. 
This contradicts unicity. Then the assumption $X_L((U_L)_{\rho_L})(L)=\{\alpha\}$ implies that $X(U_\rho)(k)=\{a\}$. Vice-versa, if $X(U_\rho)(k)=\{a\}$, then certainly 
$X(U_\rho)_L(L)$ consists of one point $\alpha\otimes L$.
\end{proof}

Observe now the following consequence of the section conjecture.
\begin{lem}\label{le:4}
Assume SC3 of Conjecture \ref{cjs:2} holds for a hyperbolic curve  $U$ and let  $V\subset U$ be an open defined over $k$.
Then any neutral fibre functor of ${\sf FC}(U)$ extends to a neutral fibre functor of ${\sf FC}(V)$.
\end{lem}
\begin{proof}
Let $\rho\in {\sf Fib}_k(U)$. Thus either (i): $\rho=\rho_u$, $u\in U(k)$ or (ii): 
$\rho=\tau_{x,\eta}$, $x\in (X\setminus U)(k)$.
In the first case, if $u\in V(k)$ then for the extension just take $\rho_u$, if $u\in (U\setminus V)(k)$
then for the extension just take any tangential fibre functor of ${\sf FC}(V)$ at $u$.
In the second case for the extension just take $\tau_{x,\eta}$ as tangential fibre functor for ${\sf FC}(V)$.
\end{proof}
This lemma suggests to pose  the following conjecture which is an immediate consequence of 
the section conjecture.
\begin{conj}\label{cje:5}
Let $U$ be a hyperbolic curve and $V\subset U$ defined over a characteristic $0$ field $k$. Then any $\rho\in {\sf Fib}_k(U)$ extends
to a functor in ${\sf Fib}_k(V)$.
\end{conj}

The reason for us to make this conjecture is the possibility of using it to reduce the section conjecture for
any curve to the case of $\P^1$ minus a $0$-dimensional subscheme defined over $k$. First we need
\begin{lem}\label{le:6}
Let $U$ be a hyperbolic curve defined over a characteristic $0$ field $k$ and $V\xrightarrow{p} U$ be a principal bundle under a finite $k$-group scheme $G$. Then $U_L$ satisfies SC3 of Conjecture \ref{cjs:2} for any finite extension $L$ of $k$
 if and only if $V_L$ does for those $L$.
\end{lem}
\begin{proof}Let $X,Y$ be the compactifications of $U,V$.
 Assume $U$ satisfies  SC3. Let $\sigma$ be a functor in ${\sf Fib}_k(V)$.
Then $\rho:=\sigma\circ p^*$ is a fibre functor in ${\sf Fib}_k(U)$. Thus $V$ appears in the pro-system
defining $U_\rho$ and hence $U_\rho\cong V_\sigma$.
This implies $X(U_\rho)=Y(V_\sigma)$.
Since $U$ satisfies the section conjecture, $X(U_\rho)$ has a unique $k$-point, thus $V$ also satisfies the
section conjecture.

 Vice-versa, assume $V_L$ satisfies SC3 for any finite extension $L$ of $k$. Let $\rho \in {\sf Fib}_k(U)$. Let $S$ denote
 the full subcategory of ${\sf FC}(U)$ generated by $p_*\sO_V$, which is finite
in the sense of \cite[Definition~2.5]{EP2}. Then the functor $\eta:=H^0_{DR}(V,p^*(-))$ is a neutral fibre functor
 for $S$. On the other hand,
functor $\rho$ restricted to $S$ yields a principal bundle $q:U_S\to U$. By the finiteness
of $S$, there exists a finite extension $L$ of $k$ such that $\eta\otimes L$ and
$\rho|_{S}\otimes L$ are isomorphic. Thus $(U_S)_L$ is isomorphic to
$V_L$ over $U_L$. Consequently $V_L$ appears in the pro-system defining
$(U_L)_{\rho_L}=U_\rho\times \Spec L$.
This means the projection $s_{\rho_L}:(U_L)_{\rho_L}\to U_L$ factors
through $s_\sigma:(U_L)_{\rho_L}\to V_L$, which is the universal covering
associated to the functor $\sigma:=H^0_{DR}((U_L)_{\rho_L},s_\sigma^*(-))\in {\sf Fib}_L(V_L)$.
 By assumption on $V_L$, $X_L((U_L)_{\rho_L})=Y_L((V_L)_\sigma)$ has an
$L$-point. Now the proof of Proposition \ref{pe:2} applies and yields a $k$-point
of $X(U_\rho)$.
\end{proof}
\begin{cor}\label{ce:7}
Let $V\xrightarrow{p} U$ be an \'etale covering of hyperbolic curves over $k$. 
Then $U_L$ satisfies SC3 of Conjecture \ref{cjs:2} for any finite extension $L\supset k$ if and only if
$V_L$ does for those $L$.
\end{cor}
\begin{proof}
Let $W\xrightarrow rV$ be an \'etale covering such that $s:W\xrightarrow r V\xrightarrow p U$ is
a Galois covering. Let $k':=H^0_\text{DR}(W,\sO_W)$. It is a finite field extension of $k$.
After base change $k\to k'$, $W$ is the union of $k'$-curves  $W_i,i=1,\ldots,d$, $d=|\text{Gal}(k'/k)|$, 
Choose one of them, say $W_1$. Thus $W_1\to U_{k'}$ and $W_1\to V_{k'}$ are Galois coverings.
According to Lemma \ref{le:6}, for any $L\supset k'$, $U_L$ satisfies SC3  if and only if
$W_{1,L}=W_1\times_{k'}\Spec L$ does, if and only if $V_L$ does. According to Proposition \ref{pe:2},
if $U_L$ satisfies SC3  then so does $U$, the same hold of course for $V$.
The claim of the corollary now follows.

\end{proof}
\begin{prop}\label{pe:8}Assume that Conjecture \ref{cje:5} holds true. 
Then SC3 of Conjecture  \ref{cjs:2} holds true for any finite field extension
of $k$  if it holds true for hyperbolic $U\subset \P^1$
over any finite field extensions of $k$.
In particular the section conjecture holds for any number field if it holds for $\P^1\setminus \{0,1,\infty\}$ over any number field.
\end{prop}
\begin{proof}
Let $U$ be a hyperbolic curve and denote by $X$ its compactification. 
By assumption on the validity of Conjecture \ref{cje:5}, Lemma \ref{le:4} shows that removing 
a $k$-subscheme of dimension 0 
makes the problem
harder. Let $f$ be a non-trivial element of $k(X)$. It defines  a ramified covering of $f: U\to \P^1$.
Replacing $U$ by a smaller Zariski open, we can assume that $f^{-1}f(U)=U$, so $f$ is an \'etale covering
of $W=f(U)\subset \P^1$. Now Corollary \ref{ce:7} tells us that the section conjecture holds true for
 $U_L$, $L\supset k$, if it holds true for $W_L$, $L\supset k$.
 
 Let now $k$ be a number field. Then Belyi's  theorem \cite{Bel} asserts that we can choose a map $f:X\xrightarrow {}\P^1$
 in such a way that $V:=f^{-1}(\P^1\setminus\{0,1,\infty\})\subset U$. Thus Corollary \ref{ce:7} again implies
 that $U$ satisfies the section conjecture if $\P^1\setminus\{0,1,\infty\}$ does.
 
\end{proof}

\bibliographystyle{plain}
\renewcommand\refname{References}

\end{document}